\documentclass{amsart}

\usepackage{array}
\usepackage{hyperref}
\usepackage{todonotes}
\usepackage{tikz}
\usepackage{xspace}
\usepackage{color}
\usepackage[english]{babel}
\usepackage[latin1]{inputenc}
\usepackage[T1]{fontenc}
\usepackage{amsmath}
\usepackage{amssymb}
\usepackage{textcomp}
\usepackage{graphicx}
\usepackage{amsthm}
\usepackage{xcolor}
\usetikzlibrary{decorations.pathmorphing}
\usetikzlibrary{matrix,trees}

\newtheorem{thm}{Theorem}[section]
\newtheorem{prop}[thm]{Proposition}
\newtheorem{lem}[thm]{Lemma}
\newtheorem{cor}[thm]{Corollary}

\theoremstyle{definition}
\newtheorem{ex}[thm]{Example}

\newtheorem{remark}[thm]{Remark}

\newcommand{\gbs}{Gr\"{o}bner bases\xspace}

\newcommand{\gb}{Gr\"{o}bner basis\xspace}

\newcommand{\Rress}{resultants\xspace}
\newcommand{\KK}{\mathbb{K} }

\newcommand{\polring}[1]{\mathbb{K}[x_1,x_2,\ldots,x_{#1}]}

\newcommand{\smallring}[1][n]{\mathbb{K}[x_{2},x_{3},\ldots,x_{#1}]}

\newcommand{\biring}{\mathbb{K}[x,y]}
\newcommand{\uniring}{\mathbb{K}[y]}

\newcommand{\R}{\mathcal{R}}

\newcommand{\Gcd}[1]{\operatorname{gcd}\left( #1\right)}
\newcommand{\ideal}[1]{\left\langle #1\right\rangle}
\newcommand{\g}{\ensuremath{{g}}}

\newcommand{\res}[3]{\ensuremath{\operatorname{res}_{#3}\left(#1,#2\right)}}

\newcommand{\variety}[1]{\ensuremath{\mathcal{V}\left(#1\right)}}

\newcommand{\ei}[1]{I_{#1}}
\newcommand{\lcm}{\ensuremath{\operatorname{lcm}}\xspace}

\begin{document}

\title[On Computing the Elimination Ideal Using Resultants]{On Computing the 
Elimination Ideal Using Resultants with Applications to Gr\"obner Bases}

\author{Matteo Gallet}
\address{
(Matteo Gallet) Johann Radon Institute for Computational and Applied Mathematics 
(RICAM),
Austrian Academy of Sciences,
Altenberger Stra\ss e 69,
4040 Linz, Austria
}

\author{Hamid Rahkooy}
\address{
(Hamid Rahkooy) Research Institute for Symbolic Computation (RISC),
Johannes Kepler University,
Altenberger Stra\ss e 69,
4040 Linz, Austria
}

\author{Zafeirakis Zafeirakopoulos}
\address{
(Zafeirakis Zafeirakopoulos) Mathematics Department, Galatasaray University, 
Ortak\"oy 34357, 
Istanbul, Turkey}

\begin{abstract}
Resultants and Gr\"obner bases are crucial tools in studying polynomial 
elimination theory. We investigate relations between the variety of the 
resultant of two polynomials and the variety of the ideal they generate. 
Then we focus on the bivariate case, in which the elimination ideal is 
principal. We study --- by means of elementary tools --- the difference 
between the multiplicity of the factors of the generator of the elimination 
ideal and the multiplicity of the factors of the resultant.
\end{abstract}

\maketitle

\section{Introduction}\label{sec:intro}
The aim of the work presented in this paper is to study elementary relations 
between resultants and elimination ideals. 
Given an ideal~$I$ in a polynomial ring with indeterminates $x_1, \dotsc, x_n$, 
we call first elimination ideal of~$I$ the intersection $I \cap \KK[x_2, 
\dotsc, x_n]$. Understanding such ideals is part of the so-called 
elimination problem, an old and central topic in polynomial algebra.

Historically the motivation for investigating such a problem comes from the 
polynomial systems solving and the desire to reduce a system in~$n$ 
variables to another one involving less variables. In this context, 
many different tools appeared, such as resultants and Gr\"obner bases. 

The problem of defining and investigating the notion of resultant has been 
considered, among others, by Sylvester, Bezout, Dixon, Macaulay and van der 
Waerden (see~\cite{clo-uag}). Gr\"obner focused on elimination ideals 
in~\cite{groebner}. A modern view of the theory of resultants was given by 
Gelfand, Kapranov and Zelevinski in~\cite{gkz}. A survey paper by Emiris and 
Mourrain \cite{emiris-mourrain} discusses determinantal representations of 
resultants and related computational questions.

In Section~\ref{sec:elim}, we focus on the case of ideals~$I$ generated by two 
polynomials. In this setting, it is natural to consider the resultant of 
the two polynomials with respect to one of the variables. 
We recall some well-known results in 
elimination theory, and provide an affine version of the result linking the 
variety of the resultant and the projection of the variety of the ideal~$I$. 
The main result, Corollary~\ref{cor:eliminationisprojection}, shows that, if 
the resultant is not identically zero, the variety of the elimination ideal and 
the projection of the variety of the ideal coincide.

In Section~\ref{sec:mult}, we examine the relation between the multiplicity of 
each factor of the resultant of two polynomials, and the multiplicity of the 
corresponding factor in the generator of the first elimination ideal. 
In \cite{lazard} Lazard gave a structure theorem for the minimal lexicographic 
\gb of a bivariate ideal generated by any number of polynomials, which reveals 
some of the factors of the generator of the elimination ideal.
We provide examples exhibiting possible behaviour of these two 
multiplicities.

\section{Elimination for two polynomials}\label{sec:elim}

In this section we are going to investigate some relations between the zero set 
of the resultant of two polynomials, the zero set of ideal they generate, 
and the zero set of the first elimination ideal of the latter. The main result 
is Corollary~\ref{cor:eliminationisprojection}. 

For an ideal~$I$ in $\polring{n}$ --- where $\KK$ is an algebraically closed 
field --- we denote by $\variety{I}$ its associated 
variety and by~$\ei{1}$ the first elimination ideal of~$I$, i.e., 
$\ei{1} = I \cap \KK[x_{2}, x_{3}, \dotsc, x_n]$. We recall two main 
results on the connection between~$I_{1}$ and~$\variety{I}$.

For $f_1,\ldots,f_m \in \polring{n}$, we write $f_i$ in the form 
\begin{equation}
\label{eq:cofactors}
  f_i=h_i (x_2,\ldots,x_n) x_1^{N_i} + \text{ terms of $x_1$-degree less than } 
N_i,
\end{equation}
for each $1 \leq i \leq m$.
Consider the projection $\pi \colon \KK^n \rightarrow \KK^{n-1}$:
\begin{equation}
\label{eq:projection} 
  \pi \bigl( (c_1, c_2,\ldots,c_n) \bigr) = \left(c_2,c_3,\ldots,c_n\right). 
\end{equation}

\begin{thm}[Elimination Theorem, see for example {\cite[Chapter~3.2, 
Theorem~2]{clo-iva}}]
\label{thm:varofelim}
Let $I_1$ be the first elimination ideal of an ideal $I \trianglelefteq 
\polring{n}$. Then
\[
  \variety{I_{1}} = \pi \bigl( \variety{I} \bigr) \cup 
  \bigl( \variety{h_1,\ldots,h_m} \cap \variety{I_{1}} \bigr). 
\]
\end{thm}
  
Although the projection of the variety of an ideal and the variety of the 
elimination ideal are in general not the same, the latter is the Zariski 
closure of the former.

\begin{thm}[Closure Theorem, see for example {\cite[Chapter~3.2, 
Theorem~3]{clo-iva}}]
\label{thm:zariski-closure}
Let $I_1$ be the first elimination ideal of an ideal $I \trianglelefteq 
\polring{n}$. Then 
\begin{itemize}
  \item $\variety{I_{1}}$ is the smallest affine variety containing 
    $\pi \bigl( \variety{I} \bigr)$, i.e., it is the Zariski closure of 
    $\pi \bigl( \variety{I} \bigr)$.
  \item If $\variety{I} \ne \emptyset$, then there is an affine variety 
    $W \subsetneq \variety{I_{1}}$ such that $\variety{I_{1}} \setminus 
    W \subseteq \pi \bigl( \variety{I} \bigr)$.
\end{itemize}
\end{thm}

\noindent Theorem~\ref{thm:zariski-closure} implies that 
$
  dim \, \pi \bigl( \variety{I} \bigr) = dim \, \variety{ I_{1} }. 
$
We mention a few notes about the possible dimension of 
$\variety{h_1, \dotsc, h_m}~\subseteq~\KK^{n-1}$. 
In general, $dim \, \variety{h_1,\ldots,h_m}$ can range from $0$ to $dim \, 
\variety{I_1}$. Hereafter we give examples for such cases.

\begin{ex}[Top Dimensional Case]
  For an example in which the dimension of 
  $\variety{h_1, \ldots, h_m}$ is the biggest possible, 
  take $I = \ideal{x_1h, h^2 } \trianglelefteq \polring{n}$, 
  where $h \in \KK[x_2,\ldots,x_n]$ with $dim \, \variety{h} = n-2$. 
  Then $I_{1} = \ideal{h^2}$ and $\variety{I_{1}} = 
  \pi \bigl( \variety{I} \bigr) = \variety{h^2} = \variety{h}$, 
  which means that $dim \,\variety{h_1,h_2} = n-2$.
\end{ex}

\begin{ex}[Zero Dimensional Case]
  If we take $I = \ideal{x_1x_3, x_2} \trianglelefteq 
  \KK[x_1,x_2,x_3]$, then $\variety{I}$ consists of two lines, the $x_1-$axis 
  and the $x_3-$axis. The projection of these lines along the $x_1$-axis 
  gives us the $x_3$-axis, which is a line, hence of dimension~$1$. 
  However $\variety{x_3,x_2} = \{(0,0)\}$ is a point, namely of 
  dimension~$0$.
\end{ex}
    
Also in the following example we see that we can have $\variety{I_{1}} = \pi 
\bigl( \variety{I} \bigr)$, independently of the dimension of 
$\variety{h_1,\ldots,h_m}$.

\begin{ex}
  Consider the ideal
  $I = \ideal{x_1h_1, x_1h_2} \trianglelefteq \polring{n}$, 
  where $h_i \in \KK[x_2,\ldots,x_n]$. 
  Then independently of what $h_1$ and $h_2$ are, we have
  $I_{1}= \{0\}$, which means that 
  $\variety{I_{1}} = \pi \bigl( \variety{I} \bigr) = \KK^{n-1}$.
\end{ex}

Also, not necessarily $\variety{h_1,\ldots ,h_m} \subseteq \variety{I_1}$ 
is true. 
Note that $\variety{h_1,\ldots,h_m}$ is not the complement of $\pi 
\bigl(\variety{I} \bigr)$, but contains the complement. 
Moreover, the dimensions of $\variety{h_1,\ldots,h_m}$ and the complement are 
independent of each other.

\medskip
As mentioned above, we will investigate the relation between the 
first elimination ideal and the resultant. We first introduce some notation 
concerning resultants. Let $f_1, f_2 \in \polring{n}$ 
be polynomials of degree $d_1$ and $d_2$ respectively. Think of them as 
elements of $\KK[x_2, \dotsc, x_n][x_1]$ and denote by $f_{i,j}$ the 
coefficient of~$x_1^j$ in~$f_i$. Recall that the resultant of $f_1$ and $f_2$ 
with respect to~$x_1$ is defined as
\[
  \res{f_1}{f_2}{x_1} = \operatorname{det} 
  \bigl( \operatorname{Syl}(f_1,f_2) \bigr),
\]
where $\operatorname{Syl}(f_1,f_2)$ is the Sylvester matrix, namely

\[\operatorname{Syl}(f_1,f_2) =
\left(
\begin{tabular}{ >{$}c<{$} >{$}c<{$} >{$}c<{$} >{$}c<{$} >{$}c<{$} >{$}c<{$} }
f_{1,{d_1}} & \cdots & \cdots  	 & f_{1,0}   &        &         \\
	    & \ddots &        	 &           & \ddots &         \\
	    &        & f_{1,d_1} & \cdots    & \cdots & f_{1,0} \\ 
f_{2,{d_2}} & \cdots & f_{2,0}   &           &        &         \\
	    & \ddots &         	 & \ddots    &        &         \\
	    &        & \ddots    &           & \ddots &         \\
	    &        &           & f_{2,d_2} & \cdots & f_{2,0} \\ 
\end{tabular}
\right)
\hspace{-2em}
\begin{tabular}{ >{$}c<{$} }
 \left.
 \begin{tabular}{ c }
 \phantom{p} \\ \phantom{p} \\ \phantom{p} \\
 \end{tabular}
 \right\} d_2    \\
 \left. 
 \begin{tabular}{ c }
 \phantom{p} \\ \phantom{p} \\ \phantom{p} \\ \phantom{p} \\ \phantom{p} \\
 \end{tabular}
  \right\} d_1   \\
\end{tabular}
\]

In the following, we consider the connection between the zero set of the 
resultant of two polynomials~$f_1$ and~$f_2$, and the projection of the 
variety of the ideal~$\left\langle f_1, f_2 \right\rangle$. In this sense, the 
situation is similar to the one of the Elimination Theorem. The homogeneous 
version of this result is an easy consequence of the basic properties of the 
resultant. Hereafter we propose an affine version of it, based upon the ideas 
of~\cite[Chapter~3.6, Proposition~3]{clo-iva}.

\begin{prop}\label{thm:var-res}
Let $I = \ideal{f_1,f_2} \in \polring{n}$ and $\R = \res{f_1}{f_2}{x_1}$. 
Recall that we denoted by~$\pi$ the projection $\KK^n \longrightarrow 
\KK^{n-1}$ defined in Equation~\eqref{eq:projection} and let $h_1$ and $h_2$ 
be as introduced in Equation~\eqref{eq:cofactors}. Then
\[
  \variety{\R} = \variety{h_1,h_2} \cup \pi \bigl( \variety{I} \bigr).
\]
\end{prop}
\begin{proof}
  We prove the following three statements:
  \begin{enumerate}
    \item First, $\variety{h_1,h_2} \subseteq \variety{\R}$. \\
      It is easy to see from the Laplace expansion of the Sylvester matrix, 
      that the greatest common divisor of $h_1$ and $h_2$ divides $\R$.
      Thus $\variety{h_1,h_2} \subseteq \variety{\R}$.
    \item Secondly, $\pi \bigl( \variety{I} \bigr) \subseteq \variety{\R}$. \\
      If $f_1,f_2\in \mathbb{K}[x_2,\ldots,x_n][x_1]$ have positive degree 
      in~$x_1$, then $\R \in \ei{1}$ (see~\cite[Chapter~3.6, 
      Proposition~1]{clo-iva}).
      Thus $\variety{\ei{1}} \subseteq \variety{\R}$.
      From the Elimination Theorem $\pi \bigl( \variety{I} \bigr) \subseteq 
      \variety{\ei{1}}$, so the claim is proved.
    \item Thirdly, $\variety{\R} \setminus \variety{h_1,h_2} \subseteq 
      \pi \bigl( \variety{I} \bigr)$. \\
      Let $c \notin \variety{h_1,h_2}$. Then we have two cases:
      \begin{itemize}
        \item Suppose $h_1(c)\ne 0$ and $h_2(c) \ne 0$. Then it follows that
        $\R(c)= \res{f_1(x_1,c)}{f_2(x_1,c)}{x_1}$. Thus
        \[ 
          \R(c) = 0 \quad \Leftrightarrow 
          \quad \res{f_1(x_1,c)}{f_2(x_1,c)}{x_1} = 0. 
        \]
      \item Suppose that $h_1(c) \ne 0$ and $h_2(c)=0$, or $h_1(c) = 0$ and 
        $h_2(c) \ne 0$.
        Without loss of generality, assume that $h_1(c) \ne 0$ and $h_2(c)=0$.
        Let $d_2 = \deg_{x_1} f_2$ and $m = \deg f_2(x_1,c)$; then $m < d_2$. 
        From \cite[Chapter~3.6, Proposition~3]{clo-iva} we have that
        \[
          \R(c) = h_1(c)^{d_2-m} \res{f_1 (x_1 , c)}{f_2(x_1 , c)}{x_1},
        \]
        and since $h_1(c) \ne 0$, 
        \[ 
          \R(c) = 0 \quad \Leftrightarrow \quad \res{f_1 (x_1, 
          c)}{f_2(x_1,c)}{x_1} = 0.
        \]
      \end{itemize}
      So in both cases $\R(c) = 0$ if and only if $\res{f_1 (x_1, 
      c)}{f_2(x_1 ,c)}{x_1} = 0$. On the other hand,
      \begin{eqnarray*}
        c \in \pi \bigl(\variety{f_1 , f_2} \bigr) & \Leftrightarrow & 
        \exists\, c_1 \in \KK : (c_1,c)\in  \variety{f_1, f_2} \\
        & \Leftrightarrow & \exists \, c_1 \in \variety{f_1(x_1,c), f_2(x_1,c)} 
        \\
        & \Leftrightarrow & \res{f_1(x_1,c)}{f_2(x_1,c)}{x_1} = 0.
      \end{eqnarray*}
      Thus $c \in \pi \bigl(\variety{I} \bigr)$ and $\variety{\R} \setminus 
      \variety{h_1,h_2} \subseteq \pi \bigl( \variety{I} \bigr) $.
    \end{enumerate}            
 The claim follows immediately from the three statements.
\end{proof}

\begin{cor}
\label{cor:eliminationisprojection}
If $f_1, f_2 \in \biring$ and $\R = \res{f_1}{f_2}{x}$ is not identically 
zero, then 
\[
  \variety{\ei{1}} = \pi \bigl( \variety{I} \bigr). 
\]
\end{cor}
\begin{proof}
By assumption, $\R$ is a non-zero univariate polynomial. If $\R$ is constant, 
the claim is easily proved. Otherwise, $\R$ vanishes on a finite set of points. 
By Proposition~\ref{thm:var-res}, also $\pi \bigl( \variety{I} \bigr)$ is a 
finite set of points. By the Closure Theorem we have that $\variety{\ei{1}}$ is 
the Zariski closure of $\pi \bigl( \variety{I} \bigr)$. Since finite sets are 
Zariski closed, we have that $\variety{\ei{1}} = \pi \bigl( \variety{I} 
\bigr)$.
\end{proof}

\section{Multiplicities}\label{sec:mult}

In Section~\ref{sec:elim} we focused on the relations between the varieties 
associated to an ideal and its elimination ideal. Here we want to deal with 
the algebraic side of the question; in particular, we are interested in the 
relations between the multiplicities of the factors of the resultant of two 
polynomials and the multiplicities of the factors of the generator of the 
elimination ideal. 

We fix an elimination order on the polynomial 
ring~$\polring{n}$ such that $x_i \prec x_1$ for all $i \in \{2, \dotsc, n\}$. 
The celebrated Elimination 
Property of \gbs asserts that if $G$ is a \gb for an ideal~$I$ with respect to 
the fixed elimination order, then $G \cap \mathbb{K}[x_{2},x_{3},\ldots,x_n]$ 
is a \gb for the elimination ideal~$\ei{1}$ with respect to the same order.

Given two polynomials $f_1$ and $f_2$, we denote by~$S_{12}$ their
S-polynomial, which is defined as follows: 
\[ 
  S_{12} := 
    \frac{\lcm \bigl( \operatorname{lm} \left(f_1\right),\operatorname{lm}
    \left(f_2\right) \bigr)}{\operatorname{lm}\left(f_1\right)} f_1 - 
    \frac{\lcm \bigl( \operatorname{lm}\left(f_1\right),
    \operatorname{lm}\left(f_2\right) 
    \bigr)}{\operatorname{lm}\left(f_2\right)} f_2, 
\]
where $\operatorname{lm}$ denotes the leading monomial.

\begin{lem}
\label{lem:spol-factor}
Let $f_1,f_2\in\polring{n}$ and suppose that $f_1 = h {f_1}'$ and
$f_2=h {f_2}'$, for some $h, f_1', f_2'$. Denote by 
$S_{12}$ the S-polynomial of $f_1$ and $f_2$ and by $S_{12}'$ the S-polynomial 
of $f_1'$ and $f_2'$. Then 
\[
  S_{12}=h S_{12}'. 
\]
\end{lem}
\begin{proof}
 The result follows from a direct computation. Let 
$\ell_i=\operatorname{lm}(f_i)$, 
$\ell_i'=\operatorname{lm}(f_i')$ and $\ell_h=\operatorname{lm}(h)$. Let 
$\ell = \lcm(\ell_1,\ell_2)$ and $\ell'=\lcm(\ell_1',\ell_2')$. Then
 \begin{align*}
   S_{12} &= \frac{\ell}{\ell_1} f_1 - \frac{\ell}{\ell_2} f_2 \\
          &= \frac{\ell}{\ell_1} h f_1' - \frac{\ell}{\ell_2} hf_2' \; = \; 
             h( \frac{\ell}{\ell_1} f_1' - \frac{\ell}{\ell_2} f_2').
 \end{align*}
 Since $\lcm(\ell_1,\ell_2)=\ell_h \lcm(\ell_1',\ell_2')$, we have that 
$\ell=\ell'\ell_h$. Therefore $\frac{\ell}{\ell_1}=\frac{\ell'}{{\ell_1}'}$ 
  and
\[
  h \left( \frac{\ell}{\ell_1} f_1' - \frac{\ell}{\ell_2} f_2'\right) \; = \; 
  h\left( \frac{\ell'}{{\ell_1}'} f_1' - \frac{\ell'}{{\ell_2}'} f_2'\right) 
  \; = \; h S_{12}'. \qedhere
\]
\end{proof}

\begin{remark}
One could use Lemma~\ref{lem:spol-factor} in Gr\"obner bases computations.
Start by computing the greatest common divisor of each pair of generators $f_i$ 
and $f_j$ at each step. Factor the greatest  common divisor of~$f_i$ and~$f_j$ 
out of~$f_i$ and~$f_j$. Then compute the S-polynomial of $f_i$ and $f_j$ and 
reduce it with respect to the other polynomials in the basis.
Finally, multiply the result of the reduction by the greatest common divisor of 
$f_i$ and $f_j$. This approach allows computations with smaller polynomials.
\end{remark}

Lemma~\ref{lem:spol-factor} also helps us proving the next proposition.

\begin{prop}
\label{prop:res-zero}
Let $I = \ideal{f_1,f_2} \trianglelefteq \polring{n}$ and 
$\R = \res{f_1}{f_2}{x_1}$. Then
 \[
   \R \equiv 0 \quad \Leftrightarrow \quad \ei{1}=\ideal{0}.
 \]
\end{prop}

\begin{proof} ~
\begin{description}
  \item[($\Leftarrow$)] Assume that $\ei{1}=\ideal{0}$. Since $\R \in \ei{1}$ 
we have $\R\equiv 0$.
  \item[($\Rightarrow$)] Assume that $\R\equiv 0$. Then either one of $f_i$ is 
zero (for which the theorem is trivial) or $f_1$ and $f_2$ have a common factor
$h$ with $\deg_{x_1}\left(h\right)>0$. Let $S$ be the normal form of $S_{12}$ 
(after reduction with respect to $f_1$ and $f_2$). If $S=0$, then $\{f_1,f_2\}$ 
is a \gb for the ideal $I$. Since $f_1,f_2\in\polring{n}\setminus \smallring$, 
none of them is in $\ei{1}$, and by the Elimination Property of \gbs we have 
$\ei{1}=\ideal{0}$. Now assume $S\neq 0$. Let $S_{12}'$,$f_1'$,$f_2'$ and $h$ 
be as in Lemma~\ref{lem:spol-factor}, and $S'$ be the reduced form of 
$S_{12}'$ with respect to $f_1'$ and $f_2'$. From Lemma~\ref{lem:spol-factor} 
and the fact that reducing $S_{12}$ by $f_1$ and $f_2$ is equivalent to 
reducing 
$S_{12}'$ by $f_1'$ and $f_2'$, we have that $S=h S'$. Therefore in the process 
of the \gb computation by Buchberger's algorithm, all of the new polynomials 
will have $h$ as a factor, and since $h\in\polring{n}\setminus\smallring$, all 
the polynomials in the \gb will belong to $\polring{n}\setminus\smallring$. By 
the Elimination Property of \gbs we have $\ei{1}=\ideal{0}$. \qedhere
\end{description}
\end{proof}

\begin{remark}
Assume that $I = \ideal{f_1,f_2} \trianglelefteq \biring$ 
and write $f_1$ and $f_2$ in the following form
\[
  f_i = t_i + h_i x^{d_i} + \sum_{j=1}^{d_i-1} h_{i_j} x^j ,
\]
where $d_i$ is the degree of~$f_i$ with respect to~$x$,
$t_i \in \uniring$ is the trailing coefficient, $h_i \in \uniring$ 
is the leading coefficient of~$f_i$ and $h_{i_{j}} \in \uniring$ are the other 
coefficients, for $i = 1,2$.
If we expand the Sylvester matrix we find the following divisibility relations:
\begin{equation*}
  \Gcd{h_1,h_2} \, | \, \R, \quad \Gcd{t_1,t_2} \, | \, \R,
\end{equation*}
and 
\begin{equation*}
  \Gcd{h_i,t_i,h_{i_1}, \ldots, h_{i_{ \left(d_k-1\right)}}} \, | \, \R, \qquad 
\mathrm{for\ } i=1,2.
\end{equation*}
\end{remark}

From now on we consider two polynomials $f_1,f_2 \in \biring$. As already 
mentioned, if $I = \ideal{f_1, f_2}$ and $g$ denotes the generator of the 
elimination ideal~$I_1 = I \cap \KK[y]$, then the resultant $\R = 
\res{f_1}{f_2}{x}$ is a multiple of~$g$. In particular, the factors of~$g$ 
are factors of~$\R$. If we are given the resultant~$\R$ and we want to 
recover~$g$, we just need to understand what are the correct multiplicities for 
the factors of~$g$.

Let $c\in\KK$ be a root of~$g$, and let $\mu$ and $\nu$  be the multiplicities 
of the factor corresponding to~$c$ in~$g$ and~$\R$ respectively. 
Clearly $\mu \leq \nu$. We 
exhibit some examples of situations that can arise. 

\begin{description}
  \item[Case $\nu = 1$]
Here, either the linear factor vanishing on~$c$ appears in~$g$ with 
multiplicity~$1$, or it does not appear at all. Using the notation of
Equation~\eqref{eq:cofactors}, such a factor appears either in~$g$ or 
in~$\Gcd{h_1,h_2}$. Hence, if $\R$ is squarefree, then $ g = \frac{\R} 
{\Gcd{h_1,h_2}} $.

The following example shows such a situation. Let $f_1 = xy-1$, $f_2 = x^2 y + 
y^2-4$. Then $\R = y (y^3 - 4y + 1)$ and $I_1=\ideal{y^3 - 4y + 1}$. The 
value~$0$ is a root of~$\R$ of multiplicity~$1$, but it is not a root of~$g$. 
The $\gcd$ of~$h_1$ and~$h_2$ is~$y$, and $g = \frac{\R}{\Gcd{h_1,h_2}}$.

\vskip 0.3em
\begin{tabular}{c c}
  \begin{minipage}[c]{0.4\textwidth}
      \includegraphics[scale=1]{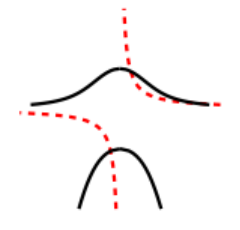}
  \end{minipage}
  &
  \begin{minipage}[c]{0.4\textwidth}
    \begin{eqnarray*}
      f_1 &=& x y - 1 \\
      f_2 &=& x^{2} y + y^{2} - 4 \\
      h_1 &=& y \\
      h_2 &=& y \\
      g   &=& y^3 - 4y + 1 \\
      \R  &=& y (y^3 - 4y + 1)
    \end{eqnarray*} 
  \end{minipage}
\end{tabular}
\vskip 0.3em

 \item[Case $\nu>1$]
Here a root of~$\R$ appears with multiplicity greater than~$1$.
\vskip 0.3em
  \begin{tabular}{c l}
    \begin{minipage}[c]{0.4\textwidth}
      \includegraphics[scale=1]{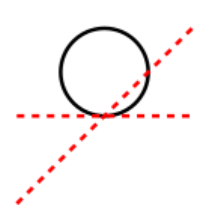}
    \end{minipage}
    &
    \begin{minipage}[c]{0.4\textwidth}
      \begin{eqnarray*}
        f_1 & = & -(y + 1) (x -  y - 1) \\
        f_2 & = & x^{2} + y^{2} - 1 \\
        h_1 & = & -(y + 1) \\
        h_2 & = & 1 \\
        g   & = & y (y + 1)^{2} \\
        \R  & = & 2 y (y + 1)^{3} \\
      \end{eqnarray*}
    \end{minipage}
  \end{tabular}
\vskip 0.3em

The factor~$y$ in~$g$ is preserved with the same multiplicity as in~$\R$, but 
the multiplicity of the factor~$(y+1)$ drops by~$1$.
\end{description}

\begin{remark}
\label{remark:drop}
Assume that no two solutions of the system given by~$f_1$ and~$f_2$ have the 
same $y$-coordinate. Suppose that the two curves defined by~$f_1$ and~$f_2$ 
admit a common tangent at an intersection point~$P$ which is parallel to the 
$x$-axis. Then the multiplicity of the factor corresponding to (the projection 
of) $P$ in~$g$ is strictly smaller than the multiplicity of the factor 
corresponding to~$P$ in~$R$.
\end{remark}

One can notice that in the case $\nu > 1$ above we are in the situation 
covered by Remark~\ref{remark:drop}, since $(y + 1)$ and the circle have a 
common tangent parallel to the $x$-axis at their intersection.

\smallskip
The multiplicity structure of isolated points can be studied 
by means of the dual space of the vanishing ideal of those points. 
In~\cite{Mantzaf-Rah-Zaf}, the problem of understanding the differences between 
the multiplicities of the factors of the resultant~$\R$ and the generator~$g$ of 
the elimination ideal has been addressed via dual spaces; there, the concept of 
directional multiplicity has been introduced to explain the exponents of the 
factors in~$\g$.

\subsection*{Acknowledgments}\label{sec:ack}

The authors would like to express their gratitude to Professors B. Buchberger, 
H. Hong, M. Kauers, and Dr. E. Tsigaridas.

The authors were supported by the strategic program "Innovatives
O\"O 2010 plus" of the Upper Austrian Government and by the Austrian 
Science Fund (FWF) grant W1214-N15, projects DK1, DK6 and DK9. Part of the 
research of the second author was carried out during a stay at the mathematics 
department of UC Berkeley supported by a Marshall Plan Scholarship. The third 
author was partially supported by Austrian Science Fund grant P22748-N18.

\bibliographystyle{plain}

\end{document}